\documentclass[a4paper,11pt]{amsart}
\usepackage[arrow,curve,matrix]{xy}

\title{Universal derived equivalences of posets of tilting modules}
\author{Sefi Ladkani}
\address{Einstein Institute of Mathematics, The Hebrew University of Jerusalem, Jerusalem 91904, Israel}
\email{sefil@math.huji.ac.il}

\DeclareMathOperator{\add}{add} \DeclareMathOperator{\rep}{rep}
\DeclareMathOperator{\Hom}{Hom} \DeclareMathOperator{\Ext}{Ext}
\DeclareMathOperator{\basic}{basic} \DeclareMathOperator{\coker}{coker}
\DeclareMathOperator{\cone}{Cone}

\newcommand{\iu}{i^{-1}}
\newcommand{\il}{i_{!}}
\newcommand{\ju}{j^{-1}}
\newcommand{\jl}{j_{*}}
\newcommand{\Qx}{Q \setminus \{x\}}
\newcommand{\vphi}{\varphi}

\newcommand{\cA}{\mathcal{A}}
\newcommand{\cD}{\mathcal{D}}
\newcommand{\cT}{\mathcal{T}}

\theoremstyle{plain}
\newtheorem*{theorem*}{Theorem}
\newtheorem{theorem}{Theorem}
\newtheorem{prop}[theorem]{Proposition}
\newtheorem{lemma}[theorem]{Lemma}
\newtheorem{cor}[theorem]{Corollary}

\numberwithin{theorem}{section}
\numberwithin{equation}{section}

\begin{document}

\begin{abstract}
We show that for two quivers without oriented cycles related by a BGP
reflection, the posets of their tilting modules are related by a simple
combinatorial construction, which we call flip-flop.

We deduce that the posets of tilting modules of derived equivalent path
algebras of quivers without oriented cycles are universally derived
equivalent.
\end{abstract}

\maketitle

\section{Introduction}

In this note we investigate the combinatorial relations between the
posets of tilting modules of derived equivalent path algebras. While it
is known that these posets are in general not isomorphic, we show that
they are related via a sequence of simple combinatorial constructions,
which we call flip-flops.

For two partially ordered sets $(X, \leq_X)$, $(Y, \leq_Y)$ and an
order preserving function $f : X \to Y$, one can define two partial
orders $\leq^f_+$ and $\leq^f_-$ on the disjoint union $X \sqcup Y$, by
keeping the original partial orders inside $X$ and $Y$ and setting
\begin{align*}
x \leq^f_{+} y \Longleftrightarrow f(x) \leq_Y y \\
y \leq^f_{-} x \Longleftrightarrow y \leq_Y f(x)
\end{align*}
with no other additional order relations. We say that two posets $Z$
and $Z'$ are related via a \emph{flip-flop} if there exist $X$, $Y$ and
$f: X \to Y$ as above such that $Z \simeq (X \sqcup Y, \leq^f_{+})$ and
$Z' \simeq (X \sqcup Y, \leq^f_{-})$.

Throughout this note, the field $k$ is fixed. Given a (finite) quiver
$Q$ without oriented cycles, consider the category of
finite-dimensional modules over the path algebra of $Q$, which is
equivalent to the category $\rep Q$ of finite dimensional
representations of $Q$ over $k$, and denote by $\cT_Q$ the poset of
tilting modules in $\rep Q$ as introduced
by~\cite{RiedtmannSchofield91}. For more information on the partial
order on tilting modules see~\cite{HappelUnger05}, the
survey~\cite{Unger07} and the references therein.

Let $x$ be a source of $Q$ and let $Q'$ be the quiver obtained from $Q$
by a BGP reflection, that is, by reverting all arrows starting at $x$.
The combinatorial relation between the posets $\cT_Q$ and $\cT_{Q'}$ is
expressed in the following theorem.

\begin{theorem} \label{t:flipflop}
The posets $\cT_Q$ and $\cT_{Q'}$ are related via a flip-flip.
\end{theorem}

In fact, the subset $Y$ in the definition of a flip-flop can be
explicitly described as the set of tilting modules containing the
simple at $x$ as direct summand, and we show that it is isomorphic as
poset to $\cT_{\Qx}$.

While two posets $Z$ and $Z'$ related via a flip-flop are in general
not isomorphic, they are \emph{universally derived equivalent} in the
following sense; for any abelian category $\cA$, the derived categories
of the categories of functors $Z \to \cA$ and $Z' \to \cA$ are
equivalent as triangulated categories, see~\cite{Ladkani07}.

For two quivers without oriented cycles $Q$ and $Q'$, we denote $Q \sim
Q'$ if $Q'$ can be obtained from $Q$ by a sequence of BGP reflections
(at sources or sinks). It is known that the path algebras of $Q$ and
$Q'$ are derived equivalent if and only if $Q \sim Q'$,
see~\cite[(I.5.7)]{Happel88}, hence by~\cite[Corollary~1.3]{Ladkani07}
we deduce the following theorem.

\begin{theorem} \label{t:unider}
Let $Q$ and $Q'$ be two quivers without oriented cycles whose path
algebras are derived equivalent. Then the posets $\cT_Q$ and $\cT_{Q'}$
are universally derived equivalent.
\end{theorem}

The paper is structured as follows. In Section~\ref{sec:source} we
study the structure of the poset $\cT_Q$ with regard to a source vertex
$x$, where the main tool is the existence of an exact functor right
adjoint to the restriction $\rep Q \to \rep (\Qx)$. For the convenience
of the reader, we record the dual statements for the case of a sink in
Section~\ref{sec:sink}. Building on these results, we analyze the
effect of a BGP reflection in Section~\ref{sec:reflect}, where a proof
of Theorem~\ref{t:flipflop} is given. We conclude by demonstrating the
theorem on a concrete example in Section~\ref{sec:example}.

\subsection*{Acknowledgement}

I would like to thank Fr\'{e}d\'{e}ric Chapoton for suggesting a
conjectural version of Theorem~\ref{t:unider} in the case of
finite-type quivers and for many helpful discussions.

\section{Tilting modules with respect to a source}
\label{sec:source}

Let $Q$ be a quiver. For a representation $M$ in $\rep Q$, denote by
$M(y)$ the vector space corresponding to a vertex $y$ and by $M(y \to
y')$ the linear transformation $M(y) \to M(y')$ corresponding to an
edge $y \to y'$ in $Q$.

Let $x$ be a source in the quiver $Q$, to be fixed throughout this
section.

\begin{lemma}
The inclusion $j : \Qx \to Q$ induces a pair $(\ju,\jl)$ of functors
\begin{align*}
\ju : \rep Q \to \rep (\Qx) && \jl : \rep (\Qx) \to \rep Q
\end{align*}
such that
\begin{equation} \label{e:jujlAdj}
\Hom_{\Qx}(\ju M, N) \simeq \Hom_{Q}(M, \jl N)
\end{equation}
for all $M \in \rep Q$, $N \in \rep(\Qx)$ (that is, $\jl$ is a right
adjoint to $\ju$).
\end{lemma}
\begin{proof}
We shall write the functors $\ju$ and $\jl$ explicitly. For $M \in \rep
Q$, define
\begin{align*}
(\ju M)(y) = M(y) && (\ju M)(y \to y') = M(y \to y')
\end{align*}
for any $y \to y'$ in $\Qx$. For $N \in \rep(\Qx)$, define
\begin{align}
\notag (\jl N)(y) &= N(y) && (\jl N)(y \to y') = N(y \to y') \\
\label{e:jlx} (\jl N)(x) &= \bigoplus_{i=1}^{m} N(y_i) && (\jl N)(x \to
y_i) = (\jl N)(x) \to N(y_i)
\end{align}
where $y_1,\dots,y_m$ are the endpoints of the arrows starting at $x$,
$(\jl N)(x) \to N(y_i)$ are the natural projections, and $y, y'$ are in
$\Qx$.

Now~\eqref{e:jujlAdj} follows since the maps $M(y_i) \to N(y_i)$ for $1
\leq i \leq m$ induce a unique map $M(x) \to N(y_1) \oplus \dots \oplus
N(y_m)$ such that the diagrams
\[
\xymatrix{
{M(x)} \ar[r] \ar[d] & {(\jl N)(x) = \bigoplus_{i=1}^{m} N(y_i)} \ar[d] \\
{M(y_i)} \ar[r] & {N(y_i)}
}
\]
commute for all $1 \leq i \leq m$.
\end{proof}

\begin{lemma}
The functor $\jl$ is fully faithful and exact.
\end{lemma}
\begin{proof}
Observe that $\ju \jl$ is the identity on $\rep(\Qx)$, hence for $N, N'
\in \rep(\Qx)$,
\[
\Hom_{Q}(\jl N, \jl N') \simeq \Hom_{\Qx}(\ju \jl N, N') =
\Hom_{\Qx}(N, N')
\]
so that $\jl$ is fully faithful. Its exactness follows
from~\eqref{e:jlx}.
\end{proof}

Denote by $\cD^b(Q)$ the bounded derived category $\cD^b(\rep Q)$. The
exact functors $\ju$ and $\jl$ induce functors
\begin{align*}
\ju : \cD^b(Q) \to \cD^b(\Qx) && \jl : \cD^b(\Qx) \to \cD^b(Q)
\end{align*}
with
\begin{equation}
\label{e:jujlAdjD} \Hom_{\cD^b(\Qx)}(\ju M, N) \simeq
\Hom_{\cD^b(Q)}(M, \jl N)
\end{equation}
for all $M \in \cD^b(Q)$, $N \in \cD^b(\Qx)$.

Let $S_x$ be the simple (injective) object of $\rep Q$ corresponding to
$x$.

\begin{lemma} \label{l:Sxperp}
The functor $\jl$ identifies $\rep (\Qx)$ with the right perpendicular
subcategory
\begin{equation} \label{e:Sxperp}
S_x^{\perp} = \left\{ M \in \rep Q \,:\, \Ext^i(S_x, M) = 0 \text{ for
all $i \geq 0$} \right\}
\end{equation}
of $\rep Q$.
\end{lemma}
\begin{proof}
Observe that $\ju S_x = 0$. Hence by~\eqref{e:jujlAdjD},
\[
\Ext^i_Q(S_x, \jl N) = \Ext^i_{\Qx}(\ju S_x, N) = 0
\]
for all $N \in \rep(\Qx)$.

Conversely, let $M$ be such that $\Ext^i_Q(S_x, M) = 0$ for $i \geq 0$,
and let $\vphi : M \to \jl \ju M$ be the adjunction morphism. From $\ju
\jl \ju M = \ju M$ we see that $(\ker \vphi)(y) = 0 = (\coker
\vphi)(y)$ for all $y \neq x$.

From $0 \to \ker \vphi \to M$ we get
\begin{equation} \label{e:kerphi0}
0 \to \Hom_Q(S_x, \ker \vphi) \to \Hom_Q(S_x, M) = 0
\end{equation}
hence $\ker \vphi = 0$. Thus $0 \to M \to \jl \ju M \to \coker \vphi
\to 0$ is exact, and from
\[
0 = \Hom_Q(S_x, \jl \ju M) \to \Hom_Q(S_x, \coker \vphi) \to
\Ext^1_Q(S_x, M) = 0
\]
we deduce that $\coker \vphi = 0$, hence $M \simeq \jl \ju M$.
\end{proof}

\begin{lemma} \label{l:jlindec}
The functor $\jl$ takes indecomposables of $\rep (\Qx)$ to
indecomposables of $\rep Q$.
\end{lemma}
\begin{proof}
Let $N$ be an indecomposable representation of $\Qx$, and assume that
$\jl N = M_1 \oplus M_2$. Then $N \simeq \ju \jl N = \ju M_1 \oplus \ju
M_2$, hence we may assume that $\ju M_2 = 0$.

Thus $M_2 = S_x^n$ for some $n \geq 0$. But $\jl N$ belongs to the
right perpendicular subcategory $S_x^{\perp}$ which is closed under
direct summands, hence $n=0$ and $M_2 = 0$.
\end{proof}

Recall that $T \in \rep Q$ is a \emph{tilting module} if
$\Ext^i(T,T)=0$ for all $i>0$, and the direct summands of $T$ generate
$\cD^b(Q)$ as a triangulated category. If $T$ is basic, the latter
condition can be replaced by the condition that the number of
indecomposable summands of $T$ equals the number of vertices of $Q$.

For a tilting module $T$, define
\[
T^{\perp} = \left\{ M \in \rep Q \,:\, \Ext^i(T, M) = 0 \text{ for all
$i > 0$} \right\}
\]
and set $T \leq T'$ if $T^{\perp} \supseteq T'^{\perp}$.
By~\cite{HappelUnger05}, $T \leq T'$ if and only if $\Ext^i_Q(T, T') =
0$ for all $i > 0$.

Denote by $\cT_Q$ the set of basic tilting modules of $\rep Q$, and by
$\cT_Q^x$ the subset of $\cT_Q$ consisting of all tilting modules which
have $S_x$ as direct summand.

\begin{lemma}
$\cT_Q^x$ is an open subset of $\cT_Q$, that is, if $T \in \cT_Q^x$ and
$T \leq T'$, then $T' \in \cT_Q^x$.
\end{lemma}
\begin{proof}
Let $T \in \cT_Q^x$ and $T' \in \cT_Q$ such that $T \leq T'$. Then $T'
\in T^{\perp}$, and in particular $\Ext^i(S_x, T') = 0$ for $i > 0$.
Since $S_x$ is injective, it follows that $\Ext^i(T', S_x) = 0$ for $i
> 0$, hence if $T' \not \in \cT_Q^x$, then $S_x \oplus T'$ would also be a
basic tilting module, contradiction to the fact that the number of
indecomposable summands of a basic tilting module equals the number of
vertices of $Q$.
\end{proof}

\begin{prop} \label{p:juTtilt}
Let $T$ be a tilting module in $\rep Q$. Then $\ju T$ is a tilting
module of $\rep (\Qx)$.
\end{prop}
\begin{proof}
We consider two cases. First, assume that $T$ contains $S_x$ as direct
summand. Write $T = S_x^n \oplus T'$ with $n > 0$, where $T'$ does not
have $S_x$ as direct summand. Then $\ju T = \ju T'$ and $T' \in
S_x^{\perp}$, hence $\jl \ju T' = T'$ and
\begin{multline} \label{e:extjuTSx}
\Ext^i_{\Qx}(\ju T, \ju T) = \Ext^i_{\Qx}(\ju T', \ju T')
\\ = \Ext^i_Q(T', \jl \ju T') = \Ext^i_Q(T', T') = 0
\end{multline}

Now assume that $T$ does not contain $S_x$ as direct summand, and let
$\vphi : T \to \jl \ju T$ be the adjunction morphism. Then $\Hom_Q(S_x,
T) = 0$ and similarly to~\eqref{e:kerphi0}, we deduce that $\ker \vphi
=0$. Observe that $\coker \vphi = S_x^n$ for some $n \geq 0$ is
injective, hence from the exact sequence $0 \to T \to \jl \ju T \to
\coker \vphi \to 0$ we get for $i > 0$,
\begin{equation} \label{e:extjuT}
0 = \Ext^i(T, T) \to \Ext^i(T, \jl \ju T) \to \Ext^i(T, \coker \vphi) =
0
\end{equation}
therefore $\Ext^i_{\Qx}(\ju T, \ju T) = \Ext^i_Q(T, \jl \ju T) = 0$ for
$i > 0$.

To show that the direct summands of $\ju T$ generate $\cD^b(\Qx)$, it
is enough to verify that for any $y \in \Qx$, the corresponding
projective $P_y$ in $\rep(\Qx)$ has a resolution with objects from
$\add \ju T$. Indeed, let $y \in \Qx$ and consider the projective
$\widetilde{P}_y$ of $\rep Q$. Applying the exact functor $\ju$ on an
$\add T$-resolution of $\widetilde{P}_y$ gives the required $\add \ju
T$-resolution of $P_y = \ju \widetilde{P}_y$.
\end{proof}

Note that $\ju T$ may not be basic even if $T$ is basic. Write
$\basic(\ju T)$ for the module obtained from $\ju T$ by deleting
duplicate direct summands. Then $\basic(\ju T)$ is a basic tilting
module with $\basic(\ju T)^{\perp} = (\ju T)^{\perp}$. It follows by
the adjunction~\eqref{e:jujlAdjD} that for $N \in \rep(\Qx)$,
\[
N \in (\ju T)^{\perp} \Longleftrightarrow \jl N \in T^{\perp}
\]

\begin{cor}
The map $\pi_x: T \mapsto \basic(\ju T)$ is an order-preserving
function $(\cT_Q, \leq) \to (\cT_{\Qx}, \leq)$.
\end{cor}
\begin{proof}
Let $T \leq T'$ and consider $N \in (\ju T')^{\perp}$. Then $\jl N \in
T'^{\perp} \subseteq T^{\perp}$, hence $N \in (\ju T)^{\perp}$, so that
$\ju T \leq \ju T'$.
\end{proof}

Let $N, N'$ be objects of $\rep (\Qx)$ with $\Ext^i_{\Qx}(N,N')=0$ for
all $i > 0$. By the adjunctions~\eqref{e:jujlAdjD},
\[
\begin{split}
& \Ext^i_Q(\jl N, \jl N') \simeq \Ext^i_{\Qx}(\ju \jl
N, N') = \Ext^i_{\Qx}(N, N') = 0 \\
& \Ext^i_Q(S_x, \jl N') \simeq \Ext^i_{\Qx}(\ju S_x, N') = 0 \\
& \Ext^i_Q(\jl N, S_x) = 0
\end{split}
\]
where the last equation follows since $S_x$ injective. Hence
\begin{equation} \label{e:extNN'}
\Ext^i_Q(S_x \oplus \jl N, S_x \oplus \jl N') = 0 \text{ for all $i>0$}
\end{equation}

\begin{cor}
Let $T$ be a basic tilting module in $\rep (\Qx)$. Then $S_x \oplus \jl
T$ is a basic tilting module in $\rep Q$.
\end{cor}
\begin{proof}
Indeed, $\Ext^i_Q(S_x \oplus \jl T, S_x \oplus \jl T) = 0$ for $i > 0$,
by~\eqref{e:extNN'}.

Let $n$ be the number of vertices of $Q$. Since $T$ is a basic tilting
module for $\Qx$, it has $n-1$ indecomposable summands, hence by
Lemmas~\ref{l:Sxperp} and~\ref{l:jlindec}, $\jl T$ decomposes into
$n-1$ indecomposable summands. It follows that $S_x \oplus \jl T$ is a
tilting module.
\end{proof}

\begin{cor}
The map $\iota_x : T \mapsto S_x \oplus \jl T$ is an order preserving
function $(\cT_{\Qx}, \leq) \to (\cT_Q^x, \leq)$.
\end{cor}
\begin{proof}
Let $T \leq T'$ in $\cT_{\Qx}$. Then $\Ext^i_{\Qx}(T,T')=0$ for all $i
>0$ and the claim follows from~\eqref{e:extNN'}.
\end{proof}

\begin{prop}
We have
\[
\pi_x \iota_x(T) = T
\]
for all $T \in \cT_{\Qx}$. In addition,
\[
T \leq \iota_x \pi_x (T)
\]
for all $T \in \cT_Q$, with equality if and only if $T \in \cT_Q^x$.
\end{prop}

In particular we see that $\iota_x$ induces a retract $\iota_x \pi_x$
of $\cT_Q$ onto $\cT_Q^x$ and an isomorphism of posets between
$\cT_{\Qx}$ and $\cT_Q^x$.

\begin{proof}
If $T \in \cT_{\Qx}$, then $\ju(S_x \oplus \jl T) = \ju \jl T = T$,
hence $\pi_x \iota_x(T) = \basic(T) = T$.

Let $T \in \cT_Q$. Then $\Ext^i_Q(T, S_x) = 0$ for $i > 0$. Moreover,
by the argument in the proof of Proposition~\ref{p:juTtilt}
(see~\eqref{e:extjuTSx} and~\eqref{e:extjuT}), $\Ext^i_Q(T, \jl \ju T)
= 0$. It follows that $S_x \oplus \jl \ju T \in T^{\perp}$, thus $T
\leq \iota_x \pi_x (T)$.

If $T = \iota_x \pi_x (T)$, then obviously $T$ has $S_x$ as summand, so
that $T \in \cT_Q^x$. Conversely, if $T \in \cT_Q^x$, then $T = S_x
\oplus T'$ with $T' \in S_x^{\perp}$, and by Lemma~\ref{l:Sxperp}, $T'
= \jl \ju T'$, hence $\iota_x \pi_x (T) = S_x \oplus \jl \ju T' = S_x
\oplus T' = T$.
\end{proof}

\begin{cor} \label{c:ffsource}
Let $X = \cT_Q \setminus \cT_Q^x$ and $Y = \cT_Q^x$. Define $f : X \to
Y$ by $f = \iota_x \pi_x$. Then $\cT_Q \simeq (X \sqcup Y,
\leq^f_{+})$.
\end{cor}
\begin{proof}
Let $T \in X$ and $T' \in Y$. If $T \leq T'$, then by the previous
proposition,
\[
f(T) = \iota_x \pi_x(T) \leq \iota_x \pi_x(T') = T'
\]
hence $T \leq T'$ in $\cT_Q$ if and only if $f(T) \leq T'$ in
$\cT_Q^x$.
\end{proof}

\section{Tilting modules with respect to a sink}
\label{sec:sink}

Now let $Q'$ be the quiver obtained from $Q$ by reflection at the
source $x$. For the convenience of the reader, we record, without
proofs, the analogous (dual) results for this case.

\begin{lemma}
The inclusion $i : \Qx \to Q'$ induces a pair $(\il,\iu)$ of functors
\begin{align*}
\iu : \rep Q' \to \rep (\Qx) && \il : \rep (\Qx) \to \rep Q'
\end{align*}
such that
\[
\Hom_{\rep(\Qx)}(N, \iu M) \simeq \Hom_{\rep Q}(\il N, M)
\]
for all $M \in \rep Q$, $N \in \rep(\Qx)$ (that is, $\il$ is a left
adjoint to $\iu$).
\end{lemma}
\begin{proof}
For $M \in \rep Q'$, define
\begin{align*}
(\iu M)(y) = M(y) && (\iu M)(y \to y') = M(y \to y')
\end{align*}
for any $y \to y'$ in $\Qx$. For $N \in \rep(\Qx)$, define
\begin{align*}
(\il N)(y) &= N(y) && (\il N)(y \to y') = N(y \to y') \\
(\il N)(x) &= \bigoplus_{l=1}^{m} N(y_l) && (\il N)(y_l \to x) = N(y_l)
\to (\il N)(x)
\end{align*}
where $y_1,\dots,y_m$ are the starting points of the arrows ending at
$x$, $N(y_l) \to (\il N)(x)$ are the natural inclusions, and $y, y'$
are in $\Qx$.
\end{proof}

\begin{lemma}
The functor $\il$ is fully faithful and exact.
\end{lemma}

Let $S'_x$ be the simple (projective) object of $\rep Q'$ corresponding
to $x$.

\begin{lemma}
The functor $\il$ identifies $\rep (\Qx)$ with the left perpendicular
subcategory
\[
^{\perp}S'_x = \left\{ M \in \rep Q' \,:\, \Ext^i(M, S'_x) = 0 \text{
for all $i \geq 0$} \right\}
\]
of $\rep Q'$.
\end{lemma}

\begin{lemma}
The functor $\il$ takes indecomposables of $\rep (\Qx)$ to
indecomposables of $\rep Q'$.
\end{lemma}

Denote by $\cT_{Q'}^x$ the subset of $\cT_{Q'}$ consisting of all
tilting modules which have $S'_x$ as direct summand.

\begin{lemma}
$\cT_{Q'}^x$ is a closed subset of $\cT_{Q'}$, that is, if $T \in
\cT_{Q'}^x$ and $T' \leq T$, then $T' \in \cT_{Q'}^x$.
\end{lemma}

\begin{prop}
Let $T$ be a tilting module in $\rep Q'$. Then $\iu T$ is a tilting
module of $\rep (\Qx)$.
\end{prop}

\begin{cor}
The map $\pi'_x: T \mapsto \basic(\iu T)$ is an order-preserving
function $(\cT_{Q'}, \leq) \to (\cT_{\Qx}, \leq)$.
\end{cor}

\begin{lemma}
Let $T$ be a basic tilting module in $\rep (\Qx)$. Then $S'_x \oplus
\il T$ is a basic tilting module of $\rep Q'$.
\end{lemma}

\begin{cor}
The map $\iota'_x : T \mapsto S'_x \oplus \il T$ is an order preserving
function $(\cT_{\Qx}, \leq) \to (\cT_{Q'}^x, \leq)$.
\end{cor}

\begin{prop}
We have
\[
\pi'_x \iota'_x(T) = T
\]
for all $T \in \cT_{\Qx}$. In addition,
\[
T \geq \iota'_x \pi'_x (T)
\]
for all $T \in \cT_{Q'}$, with equality if and only if $T \in
\cT_{Q'}^x$.
\end{prop}

\begin{cor} \label{c:ffsink}
Let $X' = \cT_{Q'} \setminus \cT_{Q'}^x$ and $Y' = \cT_{Q'}^x$. Define
$f' : X' \to Y'$ by $f' = \iota'_x \pi'_x$. Then $\cT_{Q'} \simeq (X'
\sqcup Y', \leq^{f'}_{-})$.
\end{cor}

\section{Tilting modules with respect to reflection}
\label{sec:reflect}

Let $F: \cD^b(Q) \to \cD^b(Q')$ be the BGP reflection defined by the
source $x$. For the convenience of the reader, we describe $F$
explicitly following~\cite[(IV.4, Exercise~6)]{GelfandManin03} (see
also~\cite{Ladkani07}).

Observe that a complex of representations of $Q$ can be described as a
collection of complexes $K_y$ of finite-dimensional vector spaces  for
the vertices $y$ of $Q$, together with morphisms $K_y \to K_{y'}$ for
the arrows $y \to y'$ in $Q$. Given such data, let $y_1,\dots,y_m$ be
the endpoints of the arrows of $Q$ starting at $x$, and define a
collection $\{K'_y\}$ of complexes by
\begin{align} \label{e:FK}
K'_x &= \cone \Bigl( K_x \to \bigoplus_{i=1}^{m} K_{y_i} \Bigr) \\
\notag K'_y &= K_y && y \in \Qx
\end{align}
with the morphisms $K'_y \to K'_{y'}$ identical to $K_y \to K_{y'}$ for
$y \to y'$ in $\Qx$, and the natural inclusions $K'_{y_i} = K_{y_i} \to
\cone(K_x \to \bigoplus K_{y_j}) = K'_x$ for the reversed arrows $y_i
\to x$ in $Q'$.

This definition can be naturally extended to give a functor
$\widetilde{F}$ from the category of complexes over $\rep Q$ to the
complexes over $\rep Q'$, which induces the triangulated equivalence
$F$. The action of $F$ on complexes is given, up to quasi-isomorphism,
by~\eqref{e:FK}.

\begin{lemma}[\cite{BGP73}]
$F$ induces a bijection between the indecomposables of $\rep Q$ other
than $S_x$ and the indecomposables of $\rep Q'$ other than $S'_x$.
\end{lemma}
\begin{proof}
If $M$ is an indecomposable of $\rep Q$, then $FM$ is indecomposable of
$\cD^b(Q')$ since $F$ is a triangulated equivalence.

Now let $M \neq S_x$ be an indecomposable of $\rep Q$. The map $M(x)
\to \bigoplus_{i=1}^{m} M(y_i)$ must be injective, otherwise one could
decompose $M = S_x^n \oplus N$ for some $n > 0$ and $N$.
Using~\eqref{e:FK} we see that $FM$ is quasi-isomorphic to the stalk
complex supported on degree $0$ that can be identified with $M' \in
\rep Q'$, given by
\begin{align}
\notag M'(x) &= \coker \Bigl( M(x) \to \bigoplus_{i=1}^{m} M(y_i) \Bigr) \\
\label{e:FMy} M'(y) &= M(y) && y \in \Qx
\end{align}
\end{proof}

Note also that from~\eqref{e:FK} it follows that $FS_x = S'_x[1]$.

\begin{cor} \label{c:juiuF}
$\ju T = \iu FT$ for all $T \in \cT_Q \setminus \cT_Q^x$.
\end{cor}
\begin{proof}
This follows from~\eqref{e:FMy}, since $T$ does not have $S_x$ as
summand.
\end{proof}

\begin{cor}
$F$ induces an isomorphism of posets $\rho : \cT_Q \setminus \cT_Q^x
\to \cT_{Q'} \setminus \cT_{Q'}^x$.
\end{cor}
\begin{proof}
For $T \in \cT_Q \setminus \cT_Q^x$, define $\rho(T) = FT$. Observe
that if $T$ has $n$ indecomposable summands, so does $FT$. Moreover, if
$T, T' \in \cT_Q \setminus \cT_Q^x$, then $\Ext^i_{Q'}(FT, FT') \simeq
\Ext^i_Q(T, T')$, hence $\rho(T) \in \cT_{Q'} \setminus \cT_{Q'}^x$ and
$\rho(T) \leq \rho(T')$ if $T \leq T'$.
\end{proof}

\begin{cor} \label{c:commdiag}
We have a commutative diagram
\[
\xymatrix@=1pc{
& {\cT_Q \setminus \cT_Q^x} \ar[ddr]_{\pi_x} \ar[rr]^{\rho}_{\simeq}
\ar[ddl]_{f} & &
{\cT_{Q'} \setminus \cT_{Q'}^x} \ar[ddl]^{\pi'_x} \ar[ddr]^{f'} \\ \\
{\cT_Q^x} & &
{\cT_{\Qx}} \ar[ll]^{\simeq}_{\iota_x} \ar[rr]_{\simeq}^{\iota'_x} & &
{\cT_{Q'}^x}
}
\]
\end{cor}
\begin{proof}
We have to show the commutativity of the middle triangle, that is,
$\pi_x = \pi'_x \rho$. Indeed, let $T \in \cT_Q \setminus \cT_Q^x$.
Then $\pi_x(T) = \basic(\ju T)$, $\pi'_x \rho(T) = \basic(\iu FT)$ and
the claim follows from Corollary~\ref{c:juiuF}.
\end{proof}

\begin{theorem}
The posets $\cT_Q$ and $\cT_{Q'}$ are related via a flip-flop.
\end{theorem}
\begin{proof}
Use Corollaries~\ref{c:ffsource}, \ref{c:ffsink} and~\ref{c:commdiag}.
\end{proof}

\section{Example}
\label{sec:example}

Consider the following two quivers $Q$ and $Q'$ whose underlying graph
is the Dynkin diagram $A_4$. The quiver $Q'$ is obtained from $Q$ by
reflection at the source $4$.
\begin{align*}
Q : \xymatrix{
{\bullet_1} \ar[r] & {\bullet_2} \ar[r]
& {\bullet_3} & {\bullet_4} \ar[l]
}
&&
Q' : \xymatrix{
{\bullet_1} \ar[r] & {\bullet_2} \ar[r]
& {\bullet_3} \ar[r] & {\bullet_4}
}
\end{align*}

For $1 \leq i \leq j \leq 4$, denote by $ij$ the indecomposable
representation of $Q$ (or $Q'$) supported on the vertices
$i,i+1,\dots,j$.

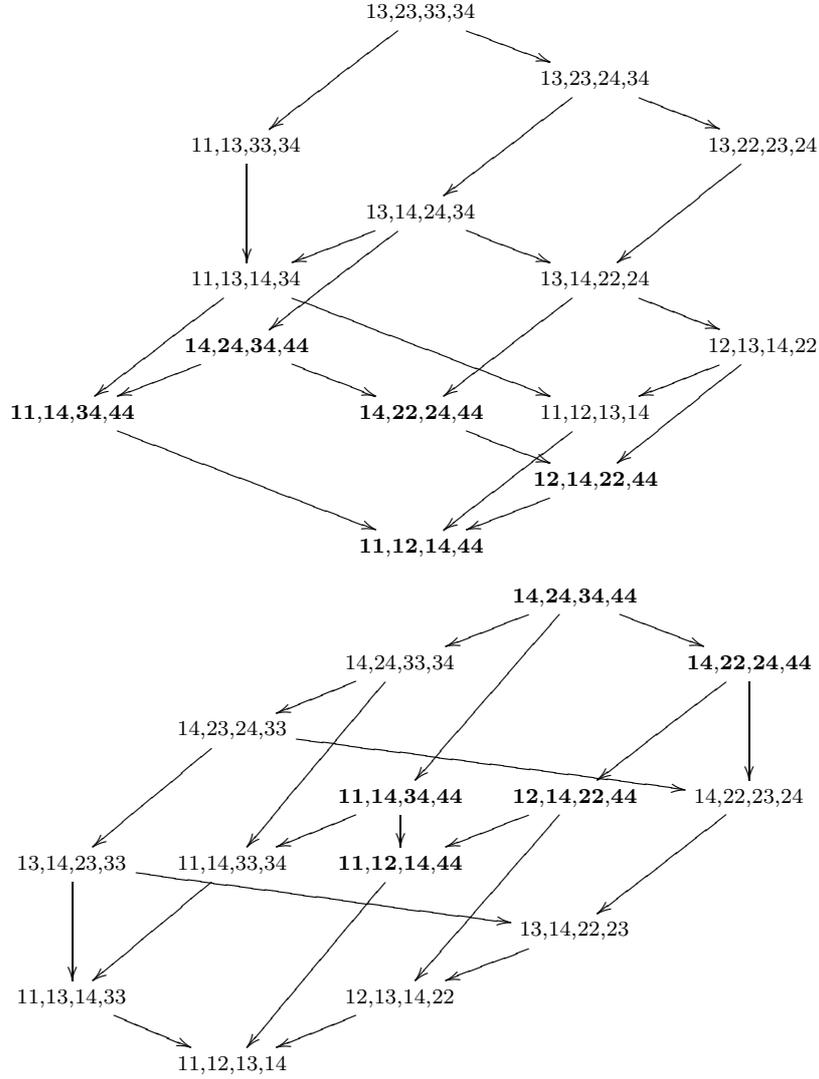
\begin{figure}
\[
\xymatrix@=1pc{
& & {_{13,23,33,34}} \ar[ddl] \ar[dr] \\
& & & {_{13,23,24,34}} \ar[ddl] \ar[dr] \\
& {_{11,13,33,34}} \ar[dd]
& & & {_{13,22,23,24}} \ar[ddl] \\
& & {_{13,14,24,34}} \ar[ddl] \ar[dl] \ar[dr] \\
& {_{11,13,14,34}} \ar[ddl] \ar[ddrr] & &
{_{13,14,22,24}} \ar[ddl] \ar[dr] \\
& {\mathbf{_{14,24,34,44}}} \ar[dl] \ar[dr] & & &
{_{12,13,14,22}} \ar[dl] \ar[ddl] \\
{\mathbf{_{11,14,34,44}}} \ar[ddrr]
& & {\mathbf{_{14,22,24,44}}} \ar[dr]
& {_{11,12,13,14}} \ar[ddl] \\
& & & {\mathbf{_{12,14,22,44}}} \ar[dl] \\
& & {\mathbf{_{11,12,14,44}}}
}
\]
\[
\xymatrix@=1pc{
& & & {\mathbf{_{14,24,34,44}}} \ar[dl] \ar[dddl] \ar[dr] \\
& & {_{14,24,33,34}} \ar[dl] \ar[dddl]
& & {\mathbf{_{14,22,24,44}}} \ar[ddl] \ar[dd] \\
& {_{14,23,24,33}} \ar[ddl] \ar[drrr] \\
& & {\mathbf{_{11,14,34,44}}} \ar[dl] \ar[d]
& {\mathbf{_{12,14,22,44}}} \ar[dl] \ar[dddl]
& {_{14,22,23,24}} \ar[ddl] \\
{_{13,14,23,33}} \ar[dd] \ar[drrr]
& {_{11,14,33,34}} \ar[ddl]
& {\mathbf{_{11,12,14,44}}} \ar[dddl] \\
& & & {_{13,14,22,23}} \ar[dl] \\
{_{11,13,14,33}} \ar[dr]
& & {_{12,13,14,22}} \ar[dl] \\
& {_{11,12,13,14}}
}
\]
\caption{Hasse diagrams of the posets $\cT_Q$ (top) and $\cT_{Q'}$
(bottom).} \label{fig:tilting}
\end{figure}

Figure~\ref{fig:tilting} shows the Hasse diagrams of the posets $\cT_Q$
and $\cT_{Q'}$, where we used bold font to indicate the tilting modules
containing the simple $44$ as summand. The subsets $\cT_Q^4$ and
$\cT_{Q'}^4$ of tilting modules containing $44$ are isomorphic to the
poset of tilting modules of the quiver $A_3$ with the linear
orientation.

Note that $\cT_Q$ was computed
in~\cite[Example~3.2]{RiedtmannSchofield91}, while $\cT_{Q'}$ is a
Tamari lattice and the underlying graph of its Hasse diagram is the
$1$-skeleton of the Stasheff associhedron of dimension 3,
see~\cite{BuanKrause04,Chapoton07}.

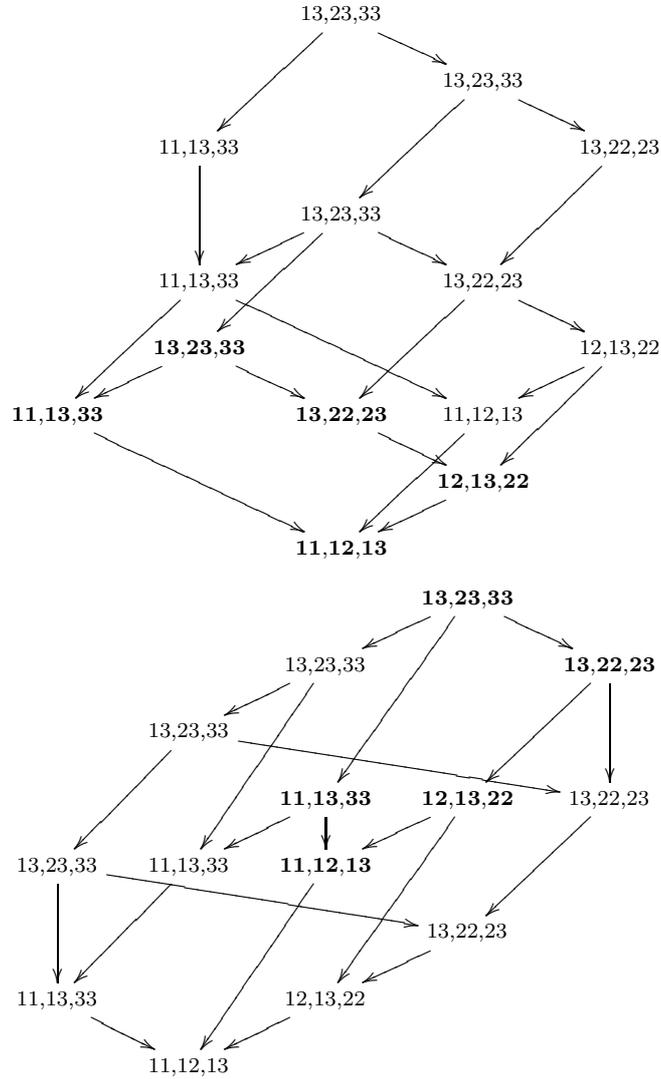
\begin{figure}
\[
\xymatrix@=1pc{
& & {_{13,23,33}} \ar[ddl] \ar[dr] \\
& & & {_{13,23,33}} \ar[ddl] \ar[dr] \\
& {_{11,13,33}} \ar[dd]
& & & {_{13,22,23}} \ar[ddl] \\
& & {_{13,23,33}} \ar[ddl] \ar[dl] \ar[dr] \\
& {_{11,13,33}} \ar[ddl] \ar[ddrr] & &
{_{13,22,23}} \ar[ddl] \ar[dr] \\
& {\mathbf{_{13,23,33}}} \ar[dl] \ar[dr] & & &
{_{12,13,22}} \ar[dl] \ar[ddl] \\
{\mathbf{_{11,13,33}}} \ar[ddrr] & &
{\mathbf{_{13,22,23}}} \ar[dr]
& {_{11,12,13}} \ar[ddl] \\
& & & {\mathbf{_{12,13,22}}} \ar[dl] \\
& & {\mathbf{_{11,12,13}}} }
\]
\[
\xymatrix@=1pc{
& & & {\mathbf{_{13,23,33}}} \ar[dl] \ar[dddl] \ar[dr] \\
& & {_{13,23,33}} \ar[dl] \ar[dddl]
& & {\mathbf{_{13,22,23}}} \ar[ddl] \ar[dd] \\
& {_{13,23,33}} \ar[ddl] \ar[drrr] \\
& & {\mathbf{_{11,13,33}}} \ar[dl] \ar[d] &
{\mathbf{_{12,13,22}}} \ar[dl] \ar[dddl]
& {_{13,22,23}} \ar[ddl] \\
{_{13,23,33}} \ar[dd] \ar[drrr] & {_{11,13,33}} \ar[ddl]
& {\mathbf{_{11,12,13}}} \ar[dddl] \\
& & & {_{13,22,23}} \ar[dl] \\
{_{11,13,33}} \ar[dr]
& & {_{12,13,22}} \ar[dl] \\
& {_{11,12,13}} }
\]
\caption{The functions $\pi_4$, $\pi'_4$ on $\cT_Q$, $\cT_{Q'}$.}
\label{fig:juiu}
\end{figure}

Figure~\ref{fig:juiu} shows the values of the functions $\pi_4$ and
$\pi'_4$ on $\cT_Q$ and $\cT_{Q'}$, respectively. The functions $f :
\cT_Q \setminus \cT_Q^4 \to \cT_Q^4$ and $f' : \cT_{Q'} \setminus
\cT_{Q'}^4 \to \cT_{Q'}^4$ can then be easily computed.

Finally, the isomorphism $\rho : \cT_Q \setminus \cT_Q^4 \to \cT_{Q'}
\setminus \cT_{Q'}^4$ is induced by the BGP reflection at the vertex
$4$, whose effect on the indecomposables (excluding $44$) is given by
\begin{align*}
11 \leftrightarrow 11 && 12 \leftrightarrow 12 && 13 \leftrightarrow 14
&& 22 \leftrightarrow 22 && 23 \leftrightarrow 24 && 33 \leftrightarrow
34
\end{align*}

\end{document}